\input amstex
\input amsppt.sty
\hoffset=-1cm \vsize=24cm \hsize=13.5cm \baselineskip=17pt
\nologo
\pageno=1
\topmatter

\def\({\bigg(}
\def\){\bigg)}
\def\Proof{\noindent{\it Proof}}
\def\Remark{\medskip\noindent{\it {Remark}}}

\title A curious identity for Bernoulli numbers\endtitle
\author
Daniel B. Gr\"unberg \qquad Hao Pan\qquad Zhi-Wei Sun
\endauthor
\abstract We prove a curious identity for the Bernoulli numbers.
\endabstract
\address
49 rue
Fondary,
75015 Paris
\endaddress
\email{daniel.b.grunberg\@gmail.com}\endemail
\address
Department of Mathematics, Nanjing University,
Nanjing 210093, People's Republic of China
\endaddress
\email{haopan1979\@gmail.com}\endemail
\address
Department of Mathematics, Nanjing University,
Nanjing 210093, People's Republic of China
\endaddress
\email{zwsun\@nju.edu.cn}\endemail
\endtopmatter
\document
\TagsOnRight

The Bernoulli numbers $B_n$ are given by
$$
\sum_{n=0}^\infty\frac{B_n}{n!}t^n=\frac{t}{e^t-1}.
$$
In [1], Gr\"unberg conjectured a curious identity for the Bernoulli numbers:
$$
\sum_{r=1}^{n-1}\frac{(-1)^rB_r}{r}\bigg(\frac{1}{n}\sum_{l=r}^{n}(-1)^l\binom{n}{l}H_{l-1}+\frac{1}{r(n-r)}\bigg)=\frac{H_{n-1}^{(2)}}{n}+\frac{H_{n-1}}{n^2},
\tag 1
$$
where
$$
H_n=\sum_{k=1}^n\frac1{k},\qquad H_n^{(2)}=\sum_{k=1}^n\frac1{k^2}.
$$
Notice that when $n=p$ is prime, in view of the Wolstenholme congruence and the von Staudt-Clausen theorem, both sides of (1) are $p$-integral. So this identity is surely reasonable. In this short note, we shall confirm this conjectured identity.
\proclaim{Theorem}
(1) holds.
\endproclaim
However, it seems not easy to deduce (1) from the generating function of $B_n$. So we shall use a difference-differential method to prove Theorem 1.
Note that
$$
\align
\sum_{l=r}^{n}(-1)^l\binom{n}{l}H_{l-1}=&\sum_{l=r}^{n}(-1)^l\bigg(\binom{n-1}{l}+\binom{n-1}{l-1}\bigg)H_{l-1}\\
=&\sum_{l=r}^{n-1}(-1)^l\binom{n-1}{l}H_{l-1}
-\sum_{l=r-1}^{n-1}(-1)^l\binom{n-1}{l}H_{l}\\
=&-\sum_{l=r}^{n-1}\frac{(-1)^l}{l}\binom{n-1}{l}
+(-1)^r\binom{n-1}{r-1}H_{r-1}.
\endalign
$$
Evidently (1) is equivalent to
$$
\align
\sum_{r=1}^{n-1}(-1)^r\frac{B_r}{r}\bigg(\frac{n}{r(n-r)}-\sum_{l=r}^{n-1}\frac{(-1)^l}{l}\binom{n-1}{l}+(-1)^r
\binom{n-1}{r-1}H_{r-1}\bigg)
=H_{n-1}^{(2)}+\frac{H_{n-1}}{n}.
\endalign
$$
Furthermore, we have
$$
\align
\sum_{l=r}^{n}\frac{(-1)^l}{l}\binom{n}{l}=&
\sum_{l=r}^{n-1}\frac{(-1)^l}{l}\binom{n-1}{l}+
\frac{1}{n}\sum_{l=r}^{n}(-1)^l\binom{n}{l}\\=&
\sum_{l=r}^{n-1}\frac{(-1)^l}{l}\binom{n-1}{l}+
\frac{(-1)^{r}}{n}\binom{n-1}{r-1}.
\endalign
$$
Hence using induction on $n$, we only need to prove that
$$
\align
&\sum_{r=1}^{n-1}\frac{B_r}{r}
\binom{n-1}{r-2}H_{r-1}-\sum_{r=1}^{n-1}\frac{(-1)^rB_r}{r(n+1-r)(n-r)}
-\sum_{r=1}^{n}\frac{B_r}{r}\frac{1}{n}\binom{n-1}{r-1}\\
=&\frac{1}{n^2}+\frac{H_n}{n+1}-\frac{H_{n-1}}{n}-(-1)^n\frac{B_n}{n^2}-(-1)^n\frac{B_{n}}{n}-B_nH_{n-1}.\tag2
\endalign
$$
With help of the recurrence relation
$$
(-1)^nB_n=\sum_{r=0}^{n}B_r\binom{n}{r},
$$
it suffices to prove that
$$\align
&\sum_{r=1}^{n-1}\frac{B_r}{r}
\binom{n-1}{r-2}H_{r-1}-
\sum_{r=1}^{n-1}\frac{(-1)^rB_r}{r(n+1-r)(n-r)}\\
=&\frac{H_n}{n+1}-\frac{H_{n-1}}{n}-(-1)^n\frac{B_{n}}{n}-H_{n-1}B_n.\tag3
\endalign$$

The  Bernoulli polynomial $B_n(x)$ are given by
$$
\sum_{n=0}^\infty\frac{B_n(x)}{n!}t^n=\frac{te^{xt}}{e^t-1}.
$$
Clearly $B_n=B_n(0)$. And it is well-known that
$$
\frac{d}{d x}B_{n+1}(x)=(n+1)B_n(x),\quad B_n(1-x)=(-1)^nB_n(x),\quad\Delta_x(B_{n+1}(x))=(n+1)x^n,
$$
where the difference operator $\Delta_x$ is defined by
$$
\Delta_x(f(x))=f(x+1)-f(x).
$$
\proclaim{Lemma}
$$\align
&\sum_{r=2}^{n+1}\frac{(-1)^rB_{r}(x)}{r}
\binom{n-1}{r-2}\binom{y}{r-1}\\
=&\sum_{r=1}^{n}
\binom{n-1}{r-1}\binom{y+n-r}{n}\bigg(\frac{B_{r+1}(-x)}{r+1}+\frac{B_{r}(-x)}{r}\bigg)+\frac{1}{y+1}\binom{y+n}{n+1}.\tag 4
\endalign$$
\endproclaim
\Proof. We have
$$\align
\Delta_x\bigg(\sum_{r=2}^{n+1}\frac{B_{r+1}(x)}{r(r+1)}
\binom{n-1}{r-2}\binom{y}{r-1}\bigg)=&\sum_{r=1}^{n+1}\frac{x^{r}}{r}
\binom{n-1}{r-2}\binom{y}{r-1}.\endalign
$$
We shall show that
$$\align
&\sum_{r=2}^{n+1}\frac{x^{r}}{r}
\binom{n-1}{r-2}\binom{y}{r-1}
\\=&\frac{1}{y+1}\binom{y+n}{n+1}+\sum_{r=1}^{n}
\binom{n-1}{r-1}\binom{y+n-r}{n}\bigg(\frac{(x-1)^{r+1}}{r+1}+\frac{(x-1)^{r}}{r}\bigg).\tag 5
\endalign$$
Note that
$$\align
&\sum_{r=0}^{n}x^{r}
\binom{n}{r}\binom{y}{r+d}=\sum_{r=0}^{n}
\binom{n}{r}\binom{y}{r+d}\sum_{j=0}^{r}\binom{r}{j}(x-1)^j\\
=&\sum_{j=0}^{n}(x-1)^{j}\binom{n}{j}\sum_{r=j}^n\binom{n-j}{r-j}\binom{y}{r+d}=
\sum_{j=0}^{n}(x-1)^{j}\binom{n}{j}\binom{y+n-j}{n+d}.\tag 6
\endalign$$
Hence we have
$$\align
&\sum_{r=2}^{n+1}\frac{1}{r}
\binom{n-1}{r-2}\binom{y}{r-1}=\frac1{y+1}\sum_{r=0}^{n-1}
1^{r}\cdot\binom{n-1}{r}\binom{y+1}{r+2}=\frac1{y+1}\binom{y+n}{n+1},
\endalign$$
i.e., (5) holds for $x=1$.
Thus by taking derivative of (5) with respect to $x$, (5) is reduced to
$$
\sum_{r=2}^{n+1}x^{r-1}
\binom{n-1}{r-2}\binom{y}{r-1}=
x\sum_{r=1}^{n}
\binom{n-1}{r-1}\binom{y+n-r}{n}(x-1)^{r-1},
$$
which is also true in view of (6).

It is easy to see that for two polynomials $f_1(x)$ and $f_2(x)$,
$$
\Delta_x(f_1(x))=\Delta_x(f_2(x))\ \Longleftrightarrow\ \frac{d}{d x}f_1(x)=\frac{d}{d x}f_2(x).
$$
So we obtain that
$$\align
&\sum_{r=2}^{n+1}\frac{B_{r}(x)}{r}
\binom{n-1}{r-2}\binom{y}{r-1}=\frac{\partial}{\partial x}\bigg(\sum_{r=2}^{n+1}\frac{B_{r+1}(x)}{r(r+1)}
\binom{n-1}{r-2}\binom{y}{r-1}\bigg)\\
=&\frac{1}{y+1}\binom{y+n}{n+1}+\sum_{r=1}^{n}
\binom{n-1}{r-1}\binom{y+n-r}{n}\frac{\partial}{\partial x}\bigg(\frac{B_{r+2}(x-1)}{(r+1)(r+2)}+\frac{B_{r+1}(x-1)}{r(r+1)}\bigg)\\
=&\frac{1}{y+1}\binom{y+n}{n+1}+\sum_{r=1}^{n}
\binom{n-1}{r-1}\binom{y+n-r}{n}\bigg(\frac{B_{r+1}(x-1)}{r+1}+\frac{B_{r}(x-1)}{r}\bigg).
\endalign$$
Replacing $x$ by $1-x$ in the above equation, we get the desired result.
\qed

Now let us return the proof of (3). It is not difficult to check that
$$
\frac{d}{dy}\binom{y}{r-1}\bigg|_{y=-1}=(-1)^{r}H_{r-1},
$$
$$
\frac{d}{dy}\bigg(\frac{1}{y+1}\binom{y+n}{n+1}\bigg)\bigg|_{y=-1}=\frac{1-H_{n-1}}{n(n+1)},
$$
and
$$
\frac{d}{dy}\binom{y+n-r}{n}\bigg|_{y=-1}=\frac{(-1)^r}{n\binom{n-1}{r}}
$$
for $1\leq r\leq n-1$.
Hence taking derivative of (4) with respect to $y$, we have
$$\align
\sum_{r=1}^{n+1}\frac{B_{r}(x)}{r}
\binom{n-1}{r-2}H_{r-1}
=&\sum_{r=1}^{n-1}(-1)^r
\frac{\binom{n-1}{r-1}}{n\binom{n-1}{r}}\bigg(\frac{B_{r+1}(-x)}{r+1}+\frac{B_{r}(-x)}{r}\bigg)\\&+
(-1)^{n+1}H_n\bigg(\frac{B_{n+1}(-x)}{n+1}+\frac{B_{n}(-x)}{n}\bigg)+
\frac{1-H_{n-1}}{n(n+1)}.
\endalign$$
Clearly
$$
\frac{\binom{n-1}{r-1}}{n\binom{n-1}{r}}-\frac{\binom{n-1}{r-2}}{n\binom{n-1}{r-1}}=\frac{1}{(n-r)(n-r+1)}.
$$
It follows that
$$\align
\sum_{r=1}^{n+1}\frac{B_{r}(x)}{r}
\binom{n-1}{r-2}H_{r-1}
=&\sum_{r=1}^{n-1}
\frac{(-1)^rB_{r}(-x)}{r(n-r)(n-r+1)}+\frac{(-1)^{n-1}(n-1)B_n(-x)}{n^2}\\
&+
(-1)^{n+1}H_n\bigg(\frac{B_{n+1}(-x)}{n+1}+\frac{B_{n}(-x)}{n}\bigg)+
\frac{1-H_{n-1}}{n(n+1)}.\tag 7
\endalign$$
Thus (3) is derived by substituting $x=0$ in (7) and noting that $B_n=0$ for odd $n\geq 3$.

\Remark. The Euler polynomial $E_n(x)$ is given by
$$
\sum_{n=0}^\infty\frac{E_n(x)}{n!}t^n=\frac{2e^{xt}}{e^t+1}.
$$
We also have
$$
\frac{d}{d x}E_{n+1}(x)=(n+1)E_n(x),\qquad E_{n}(1-x)=(-1)^nE_n(x),\qquad\Delta_x^*(E_{n}(x))=2x^n,
$$
where
$$
\Delta_x^*(f(x))=f(x+1)+f(x).
$$
Similarly, we know (cf. [2, Lemma 3.1 (ii)]) that for two polynomial $f_1(x)$ and $f_2(x)$, $\Delta_x^*(f_1)=\Delta_x^*(f_2)$ if and only if $f_1=f_2$. So using a similar discussion, we may get a analogue of (7) for Euler polynomials:
$$\align
\sum_{r=1}^{n+1}\frac{E_{r}(x)}{r}
\binom{n-1}{r-2}H_{r-1}
=&\sum_{r=1}^{n-1}
\frac{(-1)^rE_{r}(-x)}{r(n-r)(n-r+1)}+\frac{(-1)^{n-1}(n-1)E_n(-x)}{n^2}\\
&+(-1)^{n+1}H_n\bigg(\frac{E_{n+1}(-x)}{n+1}+\frac{E_{n}(-x)}{n}\bigg)+
\frac{1-H_{n-1}}{n(n+1)}.\tag 8
\endalign$$


\Refs

\ref\key 1\by D. B. Gr\"unberg\paper On asymptotics, Stirling numbers, gamma function and polylogs \jour Results Math.\vol 49\yr 2006)\pages 89-125\endref

\ref\key 2\by H. Pan and Z.-W. Sun\paper New identities involving Bernoulli and Euler polynomials
\jour J. Combin. Theory Ser. A\vol 113\yr2006\pages156-175\endref



\endRefs
\enddocument